\documentclass[11pt]{article}
\usepackage{amsfonts,latexsym}
\newcommand{\cleqn}{\setcounter{equation}{0}}
\newcommand{\clth}{\setcounter{theorem}{0}}
\newcommand {\sectionnew}[1]{\section{#1}\cleqn\clth}
\newcommand{\beq}{\begin{equation}}
\newcommand{\eeq}{\end{equation}}
\newcommand{\beqa}{\begin{eqnarray}}
\newcommand{\eeqa}{\end{eqnarray}}
\newcommand{\beaa}{\begin{eqnarray*}}
\newcommand{\ben}{\begin{eqnarray*}}
\newcommand{\eaa}{\end{eqnarray*}}
\newcommand{\een}{\end{eqnarray*}}
\newcommand{\Nset}{\hfill\nonumber}
\newcommand{\text}{\textrm}
\newcommand \nc {\newcommand}
\nc \proof {\noindent {\em{Proof.\/ }}}
\nc \qed {$\Box$\hfill}
\newtheorem{theorem}{Theorem}[section]
\newtheorem{lemma}[theorem]{Lemma}
\newtheorem{proposition}[theorem]{Proposition}
\newtheorem{corollary}[theorem]{Corollary}
\newtheorem{definition}[theorem]{Definition}
\newtheorem{example}[theorem]{Example}
\newtheorem{remark}[theorem]{Remark}
\newtheorem{conjecture}[theorem]{Conjecture}
\newtheorem{question}[theorem]{Question}
\nc \bth[1] { \begin{theorem}\label{t#1} }
\nc \ble[1] { \begin{lemma}\label{l#1} }
\nc \bpr[1] { \begin{proposition}\label{p#1} }
\nc \bco[1] { \begin{corollary}\label{c#1} }
\nc \bde[1] { \begin{definition}\label{d#1}\rm }
\nc \bex[1] { \begin{example}\label{e#1}\rm }
\nc \bre[1] { \begin{remark}\label{r#1}\rm }
\nc \bcon[1] { \begin{conjecture}\label{con#1}\rm }
\nc \bque[1] { \begin{question}\label{que#1}\rm }
\nc {\eth} { \end{theorem} }
\nc {\ele} { \end{lemma} }
\nc {\epr} { \end{proposition} }
\nc {\eco} { \end{corollary} }
\nc {\ede} { \end{definition} }
\nc {\eex} { \end{example} }
\nc {\ere} { \end{remark} }
\nc {\econ} { \end{conjecture} }
\nc {\eque} { \end{question} }
\nc \eqref[1] {{\rm{(\ref{#1})}}}
\nc \thref[1]{Theorem \ref{t#1}}
\nc \leref[1]{Lemma \ref{l#1}}
\nc \prref[1]{Proposition \ref{p#1}}
\nc \coref[1]{Corollary \ref{c#1}}
\nc \deref[1]{Definition \ref{d#1}}
\nc \exref[1]{Example \ref{e#1}}
\nc \reref[1]{Remark \ref{r#1}}
\def \W {W_{1+\infty}}
\def \WN {\W(N)}
\def \a {\alpha}
\def \b {\beta}
\def \A {{\mathcal A}}
\def \M {{\mathcal M}}
\def \L {{\mathcal L}}
\def \O {{\mathcal O}}
\def \R {{\mathcal R}}
\def \D {{\mathcal D}}

\def \Cset {{\mathbb C}}
\def \Zset {{\mathbb Z}}
\def \Zset   {{\mathbb Z}}
\def \Nset {{\mathbb N}}

\def \ord { {\mathrm{ord}} }

\def \span { {\mathrm{span}} }

\def \deg { {\mathrm{deg}} }
\def \mult { {\mathrm{mult}} }

\def \ad { {\mathrm{ad}} }
\def \Ad { {\mathrm{Ad}} }
\def \wt { {\mathrm{wt}} }
\def \psd { {\mathrm{Psd}} }

\def \Re { {\mathrm{Re}} }
\def \p { {\partial}}
\renewcommand \ker { {\mathrm{Ker}} }
\nc \Wr {Wr}
\nc \GRN { \Gr^{(N)} }
\nc \GRA[1] { \Gr_A^{(#1)} }   
\nc \GRAN { \GRA{N} }
\nc \GrA[1] { \Gr_A(#1) }\nc \GrAa { \GrA{\alpha} }
\nc \GRB[1] { \Gr_B^{(#1)} }   
\nc \GRBN { \GRB{N} }
\nc \GrB[1] { \Gr_B(#1) }
\nc \GrBb { \GrB{\beta} }
\nc \GRMB[1] { \Gr_{MB}^{(#1)} }   
\nc \GRMBN { \GRMB{N} }
\nc \GrMB[1] { \Gr_{MB}(#1) }
\nc \GrMBb { \GrMB{\beta} }

\begin{document}

\title{{\LARGE\bf{ Fuchsian bispectral operators }}}
\author{
E. ~Horozov
\thanks{E-mail: horozov@picard.ups-tlse.fr; Current address:
Laboratoire de Math\'ematiques E. Picard, CNRS UMR 5580, UFR MIG
Universit\'e Toulouse III, 31062 Toulouse C\'edex, France}
\quad
T. ~Milanov
\thanks{E-mail: fn10622@fmi.uni-sofia.bg}
\\ \hfill\\ \normalsize \textit{
Department of Mathematics and Informatics,}\\
\normalsize \textit{ Sofia University, 5 J. Bourchier Blvd.,
Sofia 1126, Bulgaria }     }
\date{}
\maketitle
\begin{abstract}
The aim of this paper is to classify the bispectral operators of any
rank with regular singular points (the
infinite point is the most important one). We characterise them in several
ways. Probably the most important result is that they are all
Darboux transformations of powers of generalised Bessel operators
(in the terminology of \cite{BHY1}).  For this reason they can
be effectively parametrised
 by the points of a certain (infinite) family of algebraic manifolds as pointed
 out in \cite{BHY1}.
\end{abstract}
\setcounter{section}{-1}
\sectionnew{Introduction}

The present paper is devoted to the characterisation and
the classification of bispectral
operators of any rank and order with only regular singularities.
 Before stating our results and placing them
properly amongst the other research we would like to give few definitions and
to recall some of the fundamental results in the area.

An ordinary differential operator $L(x,\p_x)$ is called
bispectral if it has an eigenfunction $\psi(x,z)$, depending also on
the spectral parameter $z$, which is at the same time an eigenfunction of
another differential operator $\Lambda(z,\partial_z)$
now in the spectral
parameter $z$.  In other words we look for operators $L$ , $\Lambda$ and a
function $\psi(x,z)$ satisfying equations of the form:
\beq
L\psi=f(z)\psi, \label{1.1}
\eeq
\beq   
\Lambda\psi=\theta(x)\psi. \label{1.2}
\eeq
Initially the study of bispectral operators has been stimulated by
certain
problems of computer tomography (cf. \cite{G1,G2}).
Later it turned out that the bispectral operators are
connected to several actively developing areas of mathematics
and physics - the KP-hierarchy, infinite-dimensional Lie algebras
and their representations, particle systems, automorphisms of
algebras of differential operators,
 etc. (see e.g. \cite{BHY1, BHY4, BHY5, BW, DG, W1, W2, MZ},
 as well as the papers in the proceedings volume of the conference
 in Montr\'eal \cite{BP} ). There are also indications
for eventual connections with non-commutative algebraic geometry
\cite{W3}.

In the fundamental paper \cite{DG}
Duistermaat and Gr\"unbaum raised the problem to find all bispectral
operators and completely solved it for operators $L$ of order two. The
complete
list is as follows. If we present $L$ as a Schr\"odinger operator
$$L=(\frac{d}{dx})^2+u(x),$$
the bispectral operators, apart from the obvious Airy ($u(x)=ax$) and Bessel
($u(x)=cx^{-2}$) ones, are organised into two families of potentials   $u(x)$,
which can be obtained by finitely many "rational Darboux transformations"

(1) from $u(x)=0$,

(2) from $u(x)=-(\frac{1}{4})x^{-2}$.

Thus the classification scheme prompted by the paper \cite{DG} is by the order
of the operators. G. Wilson  \cite{W1} introduced another classification scheme
- by the rank of the bispectral operator $L$.  We recall that {\it{the rank\/}}
of the operator $L$ is the dimension of the space of the joint eigenfunctions
of all operators commuting with $L$. For example all the operators of the
family (1) have rank 1, while those of the family (2) have rank 2. In the above
cited paper \cite{W1} (see also \cite{W2}) Wilson  gave a complete
description of all
bispectral operators of rank 1 (and any order). In the terminology of Darboux
transformations (see \cite{BHY1}) all bispectral operators of rank 1 are those
obtained by rational Darboux transformations on the operators with
constant coefficients, i. e. polynomials $p(\partial_x)$. Several beautiful
connections of the bispectral operators with KdV and KP-hierarchies,
algebraic curves and
Calogero-Moser particle systems have also been revealed in
\cite{DG, W1, W2} .

We will not touch upon all results in the papers \cite{DG, W1} but we would
like
to point that in both of them the classification is split into two, more or
less independent parts. First, there is an explicit construction
of families of bispectral operators of a given class (order 2 in \cite{DG};
rank 1 in \cite{W1}) The construction can be given in terms of Darboux
transformations of "canonical" operators. The second part is to give a proof
 that, if an
operator (in the corresponding class) is a bispectral one, then it belongs to
the constructed families.

In several other papers devoted to the bispectral problem (see \cite{Z, G2,
 KRo}) the
authors deal with an analog of the first part of the problem, i.e. they
construct
new families of bispectral operators. The most complete results in that
direction have been obtained in \cite{BHY1, BHY4}. To the best of our
knowledge, all known up to
now families of bispectral operators can be constructed by the methods of
the latter papers. A challenging problem is to prove that
all the bispectral operators have already been found. A natural approach would
be to divide the differential operators into suitable classes , e. g. - by
order as in \cite{DG} or by rank and to try to isolate the bispectral ones
amongst them. But having in mind
the constructions of the fundamental papers \cite{DG, W1}, with their
different and quite involved methods, the complete classification
seems to be a difficult and lengthy project. One may try to consider the
 operators with a fixed type of singularity at infinity. Obviously, then
 there arises another  difficult problem - to determine what restrictions
on the kinds of
singularities are imposed by the condition of bispectrality.

In the present paper  we consider the class of operators with regular
singularities at infinity. In fact the main results sound much stronger.
 To explain them we introduce some definitions and notations which will be
used also throughout the paper. We are going to consider operators,
 normalized as follows:
\beq
L=\sum_{k=0}^{N} V_k(x)\p_x^k, \label{1.3}
\eeq
where the coefficient at the highest derivative $V_N=1$ and the next
coefficient $V_{N-1}=0$.
Now our  assumption is that
\beq
\lim {V_j(x)}=0, \quad j=0, \ldots ,N-1 \quad when \quad
x\to \infty.
 \label{1.4}
\eeq
\noindent
(It is well known that with the above normalization
 all coefficients of $L$ are rational functions (see \cite{DG, W1}) and hence
\eqref{1.4} makes sense).

Important examples of such operators are the generalized Bessel
operators. As we are going to use them throughout
the paper we recall the definition. Introduce the notation $D=x\p_x$.

\bde{1.1}
Generalized Bessel operators $L_\b$ are the operators
\beq
L_{\b}=x^{-N}(D-\b_1)\ldots(D-\b_N),\quad (\b_1,. . . ,\b_N)\in\Cset^N,
 \label{1.5}
\eeq
\ede

\noindent
In what follows we will call the above operators by abuse of terminology (but
for simplicity) {\it Bessel operators}.

After this preparation we can formulate the result which is the core of
the present paper.

\bth{1.1}
Let $L$ be a bispectral operator \eqref{1.3} with coefficients satisfying
 \eqref{1.4}.   Then $L$ is a monomial Darboux transformation of a
Bessel operator.
\eth

The class defined by \eqref{1.3} and \eqref{1.4} includes
essentially all the bispectral operators found in \cite{DG}:
the Bessel operators and both of the families (1) and (2), the
only exception being the Airy operator. On the other hand it
includes one of the most interesting classes , found in \cite{BHY1}.
These are the operators obtained by Darboux transformations on
powers of the Bessel operators.
This class was later characterized as follows.  In \cite{BHY2} there have
been constructed highest weight modules $\M_{\b}$ with highest weight vectors
- the corresponding to \eqref{1.5} $\tau$-functions $\tau_{\b}$. Then in
\cite{BHY5} it is shown that the $\tau$-functions in the modules $\M_{\b}$ are
exactly the $\tau$-functions of the operators which are monomial Darboux
transformations.

In the course of performing the proof of \thref{1.1} we show that the
assumptions \eqref{1.3} and \eqref{1.4} for the bispectral operator $L$ impose
further restrictions on it, which justify partially the title.

\bth{1.2}
If the bispectral operator \eqref{1.3} satisfies \eqref{1.4}, then the point
$x={\infty}$ is a regular singular point.
\eth

The proof of this theorem is probably the most involved part of our
constructions (see Sect. 3). The regularity of the finite  points
 follows indirectly from \thref{1.1}.

In Sect. 4 we give another characterization of the bispectral
operators \eqref{1.3} with the restriction \eqref{1.4}.

\bth{1.3}
Any rank $r$ bispectral operator $L$ is  $\Zset_r$-invariant.
\eth

\noindent
The result is interesting and natural by itself (cf. \cite{DG, BHY1}) but in
the present paper it is also the next step in our final goal.

Finally in Sect. 6, putting together the different pieces of our construction
in the preceding sections and using the main results of \cite{BHY1, BHY5}
we obtain the following complete characterization of the Fuchsian bispectral
 operators.

\bth{1.4}
The following conditions on the operator $L$ in the form \eqref{1.3} are
equivalent:

1) L is bispectral and satisfies \eqref{1.4};

2) L is bispectral and has only regular singular points
(i.e., L is Fuchsian);

3) L is a monomial Darboux transformation of a Bessel
operator \eqref{1.5};

4) The corresponding to $L$ $\tau$-function belongs to one of the modules
$\M_{\b}$.

\eth

\noindent
In the case when the order of $L$ is two the equivalence between 1) and 3)
contains two of the most important (and difficult) theorems of \cite{DG},
concerning the families (1) and (2) above. In that sense the present paper
represents their direct generalization.

The methods which we utilize have some resemblance to the ones used in
\cite{DG}. In particular the Darboux transformations constitute one of the
main steps of our proof. But as a whole we use different ideas. First, we
work essentially with the algebraic structure of different rings of
differential or pseudo-differential operators. Essentially we do not use
 the wave function as in \cite{DG}. This we achieve by
using the bispectral involutions on pseudodifferential operators
  in Sect. 2. In the same section we observe that a bispectral
  operator $L$ (with the restrictions \eqref{1.3} and \eqref{1.4})
  satisfies a variant of the so called ``string equation'':
\beq
[L,Q]=NL^{n+1}, \label{1.6}
\eeq
where $Q$ is an operator built out of $L$. The equation \eqref{1.6}
prepares us to use certain techniques from
differential algebra in order to study the singular point of $L$  at infinity.
In particular we use the methods invented by J. Dixmier \cite{Dx} in his
studies on the Weyl algebra. Roughly speaking one associates with each
differential operator $L$ a quasi-homogeneous polynomial $p_L(X,Y)$ in such a
way that it contains the information about the ``worst'' terms of $L$ (in our
case these are the most irregular ones). See \cite{Dx} and Sect. 3 for more
details. Then in the same section the analysis of  $p_L(X,Y)$ shows that
the assumption of irregularity of the point $x=\infty$ is
incompatible with the string equation \eqref{1.6}.

The techniques from Sect. 2 is  used also in Sect. 4 to  prove
 that the rank $r$ of the operator $L$ imposes its $\Zset_r$-invariantness.
 Using it and
the fact that the infinite point is regular it is easy to perform $\Zset_r$-
invariant Darboux transformations on $L$ in order to reduce the number $n$ in
the string equation \eqref{1.6} to $0$. This automatically gives
that the operator obtained in this way is a Bessel operator.

At the end of the introduction we point out that our method
treats all ranks and orders in one scheme. We expect that some
of its components can be useful in other classification problems.

\noindent
{\flushleft{\bf{Aknowledgements. }}} We are grateful to L. Haine and G. Wilson
 for a number of
discussions and  suggestions that lead to considerable improvements of the
paper. The first mentioned author acknowledges the kind hospitality of
Universit\'e de Louvain, Universit\'e de Nantes and
Universit\'e Toulouse III, where the last part of the  paper
was done.


\sectionnew{Preliminaries}


In this section we have collected some terminology, notations and results
relevant for the study of bispectral operators.
Our main concern is to introduce unique notation which will be used throughout
the paper and to make the paper self contained. There are also few results
which cannot be found formally elsewhere, but in fact are reformulations (in
a suitable for the present paper form) of statements from other sources.

\subsection{}


In this subsection we recall some definitions, facts and notation from Sato's
theory of KP-hierarchy \cite{S, DJKM, SW} needed in the paper.
For a complete
presentation of the theory we recommend also \cite{Di, vM}.
We start with the notion of {\it the wave operator} $K(x,\p_x)$.
This is a
pseudo-differential operator
\beq
K(x,\p_x)=1+\sum_{j=1}^{\infty}a_j(x)\p_x^{-j}, \label{2.1}
\eeq
with coefficients $a_j(x)$ which could be convergent or formal
power (Laurent) series.
In the present paper we will consider $a_j$ most often as formal
Laurent
series in $x^{-1}$. The wave operator defines the (stationary)
Baker-Akhiezer function $\psi(x,z)$:
\beq
\psi(x,z)=K(x,\p_x)e^{xz}. \label{2.2}
\eeq
>From \eqref{2.1} and \eqref{2.2} it follows that $\psi$ has the following
asymptotic expansion:
\beq
\psi(x,z)= e^{xz}(1+\sum_1^{\infty}a_j(x)z^{-j}),\quad z\to {\infty}.
                                                                \label{2.3}
\eeq
Introduce also the pseudo-differential operator $P$:
\beq
P(x,\p_x)= K\p_xK^{-1}.                                         \label{2.4}
\eeq
The following spectral property of $P$, crucial in the theory
of KP-hierarchy,
is also very important for the bispectral problem:
\beq
P\psi(x,z)=z\psi(x,z).                                          \label{2.5}
\eeq
When it happens that some power of $P$, say $P^N$, is a differential
operator,
we get that $\psi(x,z)$ is an eigenfunction of an ordinary
differential operator $L=P^N$:
\beq
L\psi=z^N{\psi}. \label{2.6}
\eeq
It is possible to introduce the above objects in many different
ways, starting
with any of them (and with other, not introduced above). For us
it would be
important also to start with given {\it differential operator} $L$:
\beq
L(x,\p_x)= \p_x^N + V_{N-2}(x)\p^{N-2} + \ldots + V_0(x). \label{2.7}
\eeq
 One can define the
pseudofferental operator $P$ as an $N$-th root of the operator $L$:
\beq
P = L^{\frac{1}{N}} = \p+ \ldots, \label{2.8}
\eeq
and the wave operator $K$ as:
\beq
LK=K {\p^N} . \label{2.9}
\eeq
An important notion, connected to an operator $L$ is the algebra
$\A_L$ of operators commuting with $L$ (see \cite{Kr, BC}). This algebra
is commutative one. The wave function $\psi(x,z)$ (defined in
\eqref{2.2}) is
a common wave function for all operators $M$ from $\A_L$:
\beq
M\psi(x,z)=g_M(z)\psi(x,z). \label{2.10}
\eeq
We define also the algebra $A_L$ of all functions $g_M(z)$ for which
\eqref{2.10} holds for some $M\in \A_L$. Obviuosly the algebras $A_L$ and
$\A_L$ are isomorphic.

Following \cite{Kr} we introduce {\it the rank of the algebra} $\A_L$ as
the greatest common divisor of the orders of the operators in $\A_L$.


\subsection{}


Here we shall briefly recall the definition of Bessel wave function and
of monomial Darboux transformations from it. For more details see
\cite{BHY1}.
Let $\beta \in \Cset^N$ be such that
\beq
\sum_{i=1}^{N}\beta_i = \frac{N(N-1)}{2}.
\label{2.15}
\eeq

\bde{bess} \cite{F, Z, BHY1}
{\it {Bessel wave function}} is called the unique wave function
$\Psi_\beta(x,z)$ depending only on $xz$ and satisfying
\beq
L_\beta (x, \partial_x) \Psi_\beta(x,z) =   z^N \Psi_\beta(x,z),
\label{2.17'2}
\eeq
where the Bessel operator $ L_{\b}(x,\p_x)$ is given by \eqref{1.5}.
\ede

\noindent
Because the Bessel wave function depends only on $xz$, \eqref{2.17'2}
implies
\beqa
&&D_x \Psi_\beta(x,z)= D_z \Psi_\beta(x,z),
\label{2.17'1}\\
&&L_\beta (z, \partial_z) \Psi_\beta(x,z) =
   x^N \Psi_\beta(x,z).
\label{2.17'3}
\eeqa

To introduce the monomial Darboux transformations of Bessel
wave functions we first recall the definition of polynomial
Darboux transformations given in \cite{BHY1}.

\bde{dt}
We say that the wave
function $\Psi$ is a {\em{Darboux transformation\/}} of the
Bessel  wave function $\Psi_\beta(x,z)$ iff
there
exist polynomials $f(z)$, $g(z)$ and differential operators
$P(x,\partial_x)$,
$Q(x,\partial_x)$ such that
\beqa
&&\Psi=\frac{1}{g(z)} P(x,\partial_x) \Psi_\beta(x,z),
\label{2.81} \\
&&\Psi_\beta(x,z)=\frac{1}{f(z)} Q(x,\partial_x) \Psi.
\label{2.82}
\eeqa
The Darboux transformation is called {\em{polynomial\/}} iff
the operator $P(x, \partial_x)$ from \eqref{2.81} has the form
\beq
P(x,\partial_x)=x^{-n}\sum_{k=0}^n p_k(x^N) D_x^k,
\label{2.84}
\eeq
where $p_k$ are rational functions, $p_n\equiv 1$.
\ede

\noindent
We will need the following two definitions of monomial Darboux
transformations.
Their equivalence is proved in \cite{BHY1}.

\bde{mon1}
We say that the wave function $\Psi(x,z)$
is a {\em monomial Darboux transformation\/} of the Bessel
wave function              %
$\Psi_\beta(x,z)$ iff it is a polynomial
Darboux
transformation of $\Psi_\beta(x,z)$ with $g(z) f(z) = z^{d N}$,
$d\in\Nset$.
Further the differential operator
\ben
L=\p^M+V_{M-2}\p^{M-2}+\ldots+V_0
\een
is a monomial Darboux transformation of $L_\b $
if the wave function corresponding to $L$ is a monomial
Darboux transformation of the wave function
corresponding to $L_\b$.
\ede

\noindent
\bde{mon2}
The wave function $\Psi(x,z)$ is a
{\em monomial Darboux transformation\/} of the Bessel wave
function
$\Psi_\beta(x,z)$ iff \eqref{2.81} holds
with
$g(z)=z^n$, $n= \ord P$ and the kernel of the operator
$P(x, \partial_x)$  has
a basis consisting of several groups of the form
\beq
\partial_y^l \Big(\sum_{k=0}^{k_0}
                  \sum_{j=0}^{\mult(\beta_i+kN)-1}
                   b_{kj} x^{\beta_i + kN} y^j
              \Big)
              \Big|_{y = \ln x},
              \quad 0 \leq l \leq j_0,           \label{mon}
\eeq
where $\mult(\beta_i + kN):=$ multiplicity of
$\beta_i+ kN$ in $\bigcup_{j=1}^N \{\beta_j+ N\Zset_{\geq0}\}$ and
$j_0 =\max\{j | b_{kj}\not=0 {\textrm{ for some }} k\}$.
\ede

\noindent
>From \deref{dt} and \deref{mon1} one immediately obtains the
following description of monomial Darboux transformations:

\ble{2.4}
The differential operator $L$ is a monomial Darboux transformation of
the Bessel operator $L_{\b}$ iff there are differential
operators $P=P(x,\p_x)$ , $Q=Q(x,\p_x)$ and numbers
$d$ , $d'$ such that
\beq
Q(x,\partial_x) P(x,\partial_x) =
L_\beta(x,\partial_x)^d,                          \label{c1}
\eeq
\beq
P(x,\partial_x) Q(x,\partial_x) =
L(x,\partial_x)^{d'},                             \label{c2}
\eeq
where the operator $P$ satisfies \eqref{2.84}.
\ele

\noindent
We will also  reformulate some results from \cite{BHY1}.
In \cite{BHY1} one can find a proof of the following statement.

\ble{77.1}
If $L_{\b}$ is a Bessel operator of order $N$ and rank $r$,
there exists a
Bessel operator $L_{\b'}$ of order $r$ such that $L_{\b}$ is
a monomial
Darboux transformation of $L_{\b'}$.
\ele

\noindent
For the proof of this lemma see the proof of Proposition 2.4
from \cite{BHY1}
(although the statement there is formulated in a  different
way). We end this subsection by reformulating (in a weaker form)
the main result, which we need from \cite{BHY1}.

\bth{2.11}
The monomial Darboux transformations of the Bessel operators
are bispectral operators.
\eth


\subsection{}

Here we recall several simple properties of bispectral
operators following
\cite{DG, W1}. As we have already mentioned in the
introduction we are going to study ordinary differential
operators $L$ of arbitrary order $N$ which are
normalised as in \eqref{1.3}, i.e. with $V_N=1$ and
$V_{N-1}=0$. Assuming that $L$ is bispectral means that we have
also another operator $\Lambda$, a wave
function $\psi(x,z)$ and two other functions $f(z)$ and
$\theta(x)$, such that the equations \eqref{1.1} and \eqref{1.2}
hold. The following lemma, due to
\cite{DG}, has been fundamental for all studies of bispectral
operators.

\ble{2.31}
There exists a number $m$, such that
\beq
(\ad \, L)^{m+1}{\theta}=0. \label{2.31}
\eeq
\ele

\noindent
For its simple proof, see \cite{DG, W1}. We will
consider that $m$ is the minimal
number with this property. An important corollary of
the above lemma is the following result.

\ble{2.32}
The functions $f(z)$ and $\theta(x)$ are polynomials.
\ele

\noindent
The next result is also contained in \cite{DG, W1}, but it is
not formulated as
a separate statement. We give its short proof following \cite{W1}.

\ble{2.33}
The coefficients $\a_j$ in the expansion \eqref{2.1} of the
wave operator $K$ are rational functions.
\ele

\proof
>From the equation \eqref{2.31} it follows that
$$(\ad \, \p_x^N)^{m+1}(K^{-1}{\theta}K)=0.$$
On the other hand the kernel of the operator
$(\ad \, \p_x^N)^{m+1}$
consists of
all pseudo-differential operators whose coefficients are polynomials
in $x$ of degree at most $m$. This gives that
\beq
\theta K=K\Theta, \label{2.32}
\eeq
with a pseudo-differential operator $\Theta$ :
\beq
\Theta= \Theta_0 +\sum_1^{\infty}\Theta_j\p_x^{-j} \label{2.33}
\eeq
whose coefficients $\Theta_j$ are polynomials of degree at most $m$.
 We have $\theta=\Theta_0$. Comparing the coefficients at
$\partial_x^{-j}$ we find
that all the coefficients $\a_j(x)$ of $K$ are rational functions.
\qed

\bre{r1}
We notice that at least one of the coefficients of $\Theta_j$
has degree {\it exactly m}, where $m$ from \leref{2.31} is minimal.
This fact will be used later.
\ere

\noindent
The last lemma has as an obvious consequence one of the few general
results,
important in all studies of bispectral operators. Noticing
that the
coefficients of $L$ are polynomials in the derivatives of
$\a_j(x)$ we get

\ble{2.34}
     The coefficients of $L$ are rational functions.
\ele


\sectionnew{Bispectral involutions and the string equation}


The condition \eqref{1.4}  for vanishing of the coefficients $V_j(x)$ of
a bispectral operator $L$ implies further restrictions on all
objects connected
to $L$ - the wave function $\psi(x,z)$, the wave operator K and the
coefficients of $L$ itself. This gives us the opportunity to define two
anti-isomorphisms $b$ and $b_1$ (``bispectral involutions'') between
the algebras of of pseudo-differential operators with coefficients -
formal Laurent
series in the variables $x^{-1}$ and $z^{-1}$. In its turn using these
anti-isomorphisms will allow us to continue our further constructions
in the rest of the paper by purely algebraic analysis on
the differential or pseudo-differential operators, avoiding the wave
function.

\subsection{Bispectral involutions}

 In the next lemma, following \cite{DG} we find the simplest
restrictions on the coefficients of the wave operator $K$ and on $L$.

\ble{3.1}
(i) The coefficients $V_j(x),j=N-2,. . . 0$ of $L$ vanish at
$\infty$ at least as $x^{-2}$.

(ii) The coefficients $\a_j, \quad j=1, \dots$ of the wave operator $K$ vanish
at least as $x^{-1}$.
\ele

\proof We are going to prove both statements simultaneously. We use
the formula
$$LK=K{\p}^N, $$
\leref{2.33} and
\leref{2.34}. Comparing the coefficients at $\p^{N-2}$ at the both sides of
the above identity we get:
$$
V_{N-2}+N\a_1^{'}=0.
$$
Having in mind  that $V_{N-2}$ is equal to the derivative of the rational
function $\a_1$ and that it vanishes at $\infty$ we see that it vanishes
at least as $x^{-2}$. Continuing in the same manner we find
$$
V_{N-3}+V_{N-2}\a_1 + \frac{N(N-1)}{2}\a_1^{''}+N\a_2^{'}=0.
$$
We see that $\a_2^{'}$ is vanishing (at least as $x^{-2}$)
and  that $V_{N-3}$ vanishes
again at least as $x^{-2}$, being a sum of such terms.
By induction we get that
 $\a_{s-1}^{'}$, $s=1,. . . N-1$ vanishes at least as $x^{-2}$ and
the same holds
for $V_{N-s}, s=2,. . ,N$ as it is a sum of products
$V_j\a_m^{(k)}$, where $N-1>j>s$,
$m=1,. . . ,N-1$ (here $\a^{(k)}$ denotes k-th derivative) and also pure
derivatives of $\a_m$.  Arguing as above we get the statement of the lemma.
\qed

Following \cite{BHY4} we will introduce an anti-isomorphism $b$
between   the
algebra $\cal{B}$  of pseudo-differential operators $P(x,\p_x)$
in the variable $x$ and the algebra $\cal{B}^{'}$
of pseudo-differential operators $R(z,\p_z)$ in the variable $z$.
More precisely $\cal B$ consists of those pseudo-differential
operators
\ben
P=\sum_k^{\infty}p_j(x^{-1})\p_x^{-j},
\een
for which there is a number $n\in \Zset$ (depending on $P$)
such that $x^np_j(x^{-1}) , j=k,k+1,\ldots$
are formal power series in $x^{-1}$.
The involution
\ben
b :{\cal B} \longrightarrow {\cal B'}
\een
is defined by
\beq
b(P)e^{xz}=Pe^{xz}=\sum_k^{\infty}
                   z^{-j}p_j(\p_z^{-1})e^{xz} ,
\quad
{\textrm{for }}
\quad P\in \cal B                                    \label{3.0}
\eeq
i. e. $b$ is just a continuation of the standard
anti-isomorphism between two copies of the Weyl algebra.
In what follows we will use also the anti-isomorphism
\beq
b_1:{\cal B} \longrightarrow {\cal B'},
\qquad    b_1(P)=b(\Ad _KP).                              \label{3.3}
\eeq
Obviously $b$ and $b_1$ can be considered as involutions of $\cal B$
and without any ambiguity we can denote the inverse isomorphisms
$b^{-1},b_1^{-1}:{\cal B'}\longrightarrow \cal B$
 by the same letters.

\bre{inv}
If we use relations \eqref{1.1} and \eqref{1.2}
to define an involution $b_1$
on the subalgebra of $\cal{B}$ generated by $L$
and $\theta$, then we have
\ben
  b_1(L)=b(K^{-1}LK)=b(\Ad _KL),
\een
\ben
  b_1(\theta)=b(K^{-1}{\theta}K)=b(\Ad _K\theta).
\een
This prompts the definition \eqref{3.3}.
\ere

\noindent
Since the operators $K$ and $\Theta$ are from
$\cal B$ we can define two operators $S$ and
$\Lambda $ as follows:
\beqa
      S(z,\p_z) & = & b(K(x,\p_x)),                      \label{3.1}\\
\Lambda(z,\p_z) & = & b(\Theta).                         \label{3.5}
\eeqa
Explicitely one has
\beq
S=\sum_{j=0}^{\infty}z^{-j}\a_j(\p_z)
 =\sum_{j=0}^{\infty}a_j(z^{-1})\p_z^{-j},
  \qquad a_0=1                                      \label{3.2}
\eeq
and also
\beq
\Lambda(z,\p_z)=\sum_{j=0}^{\infty}z^{-j}\Theta_j(\p_z)
               =\sum_{i=0}^{m}\Lambda_i(z^{-1})\p_z^i,
                                        \label{3.6}
\eeq
where $\Lambda_m \neq 0$ (see \reref{r1}) and the coefficients
$\Lambda_i$ and $a_j$ should be viewed as
formal power series. We are going to see that they are polynomials
in $z^{-1}$.

\ble{3.2}
The coefficients $a_j$ of the operator $S$ are polynomials in $z^{-1}$.
\ele

\proof
Using that
\ben
LK=K\p^N,
\een
we can apply the involution $b$ and to derive:
\ben
Sb(L)=z^N S.
\een
Rewrite in details the last formula:
\ben
(\sum_0^{\infty}a_j(z^{-1})\p_z^{-j})
               (z^N+ z^{N-2}V_{N-2}(\p_z)+. . . )=
z^N(\sum_0^{\infty}a_j(z^{-1})\p_z^{-j}) .
\een
Comparing the coefficients at $\p_z^{-j}$ for $j=2,3. . . $
and having in mind
that according to \leref{3.1}:
\ben
  V_k(\p_z)= \sum_2^{\infty} V_{k,s}\p_z^{-s},
             \quad k=0,. . . ,N-2,
\een
we obtain  relations for $a_1$ and $a_2$ in the form:
\ben
-Nz^{N-1}a_1 + \sum_0^{N-2}z^kV_{k,2} & = & 0,        \\
-2Nz^{N-1}a_2+(\sum_0^{N-2}z^kV_{k,2} + N(N-1)z^{N-2})a_1
             + \sum_0^{N-2}z^kV_{k,3} & = & 0 .
\een
We see that $a_1, a_2$ are polynomials in $z^{-1}$. By
induction we get that
any $a_s$ satisfies an equation of the form:
\ben
-sNz^{N-1}a_s+\sum_0^{N-2}z^kq_{k,s}(z^{-1})=0,
\een
where $q_{k,s}$ are already polynomials in $z^{-1}$. This
proves the lemma.
\qed

Now we are ready to show that the operator $\Lambda$ has
coefficients $\Lambda_j$, which are polynomials in $z^{-1}$.
Denote temporarily by $r$ the degree of the
polynomial $\theta$, i. e. if $\theta(x)=\theta_rz^r+\ldots$,
then $\theta_r\not = 0$.

\ble{3.3}
The coefficients $\Lambda_i$ of the operator $\Lambda$ are
polynomials in $z^{-1}$. The degree of $\theta$  $r=m$ and
\beq
  \Lambda_m=\theta_m,
\qquad \Lambda_{m-1}=\theta_{m-1}, \label{3.8}
\eeq
\ele

\proof
Using the definition \eqref{3.5} and applying the involution $b$
to the relation
$\theta(x)K=K\Theta$  we get:
\ben
  \Lambda S = S{\theta(\p_z)}
\een
As the coefficients of $\Lambda$ are expressed as differential
polynomials of
the coefficients $a_j$ of $S$ we get that $\Lambda_j$ are also
polynomials in
$z^{-1}$. Comparing the first two coefficients of the above equality
we get also \eqref{3.8} .
\qed


\subsection{The string equation}

In this subsection we are going to show that for the
bispectral operator $L$ there exists another operator
$Q$, for which the string equation \eqref{1.6}
holds.  This equation as well as other properties of
the operator $Q$
(with appropriate normalisation) would be crucial for our constructions.

In what follows we would assume that the number $m$ is divisible by $N$.
This is not a restriction since we can always replace $\Lambda$ by
$\Lambda^N$. We put $m=Nl$.

\ble{4.1}
There is a natural number $n$ such that:
\beq
Q=K^{-1}x{\p_x}^{nN+1}K, \label{4.4}
\eeq
is a differential operator. The operator
$Q$ is a solution to the string equation \eqref{1.6}.
\ele

\proof
Using the bispectral property one can write
$$
(\ad \, L)^{m-1}{\theta}=(-1)^{m-1}b_1((\ad \, z^N)^{m-1}{\Lambda}) .
$$
\noindent
Each application of the operator $\ad \, z^N$ to any
differential operator $P$ reduces its order by 1.Using the fact that
the operator
$$\Lambda={\Lambda_m}\p_z^m+\Lambda_{m-2}\p_z^{m-2 } + . . . , $$
where $\Lambda_m$ is a nonzero constant, we get that the operator
$$
(\ad \, z^N)^{m-1}{\Lambda}={\Lambda_m}(\ad \, z^N)^{m-1}{\p_z^m} .
$$
is an operator of order 1.
Now prescribing weights to $z$ and to $\p_z$ as follows: $\wt (z)=1$,
$\wt (\p_z)=-1$ we obtain that the right-hand side of the above
identity has weight equal to $(m-1)N-m$. This
shows that the operator in the above equality has the form:
$$
(\ad \, z^N)^{m-1}{\Lambda}=cz^{(m-1)(N-1)}\p_z+c_1{z^{mN-m-N}},
c\not =0
$$
In this way we get that
$$
Q_1:=(\ad \, L)^{m-1}{\theta}=
b_1\Big(
        (-1)^{m-1}(cz^{(m-1)(N-1)}\p_z+c_1{z^{mN-m-N}})
   \Big)
$$
is a differential operator.
Using the fact that $m=Nl$ and that $b_1(z)=L^\frac{1}{N}$
we obtain
$$
((-1)^{m-1}Q_1-c_1L^{Nl-l-1})=
cb_1(z^{(m-1)(N-1)}\p_z)=cb_1(z^{nN+1}\p_z),
$$
where we have put $n=l(N-1)-1$. Now it is obvious that
$$
Q:=b_1(z^{nN+1}\p_z)=
\frac{1}{c}((-1)^{m-1}Q_1-c_1L^{Nl-l-1})
$$
is a differential operator.
The identity
\eqref{1.6} is obtained by applying the bispectral involution
to
\ben
[z^{nN+1}\p_z,z^N]=Nz^{N(n+1)}.
\een
\qed

\bco{3.1}
For any positive integer $i$ the following formula holds:
\beq
    (\ad \, L)^{i}(Q^i)=i! N^i L^{i(n+1)}.            \label{4.5}
\eeq
\eco

\proof
Assume that \eqref{4.5} is true for $1,2,\ldots,i$. Then
\ben
(\ad \, L)^{i+1}(Q^i)=0
\een
Since $\ad \, L$ is a differentiation in the ring of differential
operators with rational coefficients we can use the Leibnitz's rule:
\ben
(\ad \, L)^{i+1}(Q^{i+1})=(\ad \, L)^{i+1}(Q^i. Q)=\sum_{j=0}^{i+1}
{{i+1}\choose j}(\ad \, L)^{i+1-j}(Q^i)(\ad \, L)^{j}(Q)
\een
The only nonzero term in the above sum is the one for
$j=1$, hence
\ben
(\ad \, L)^{i+1}(Q^{i+1})=(i+1)N. (\ad \, L)^{i}(Q^i)
L^{n+1}.
\een
Now  \eqref{4.5} follows by induction on $i$.
\qed


\sectionnew{The infinite point}

The present section is divided into three subsections, corresponding
to the basic results, which we shortly describe. In the first
subsection we present the operator $Q$ as a polynomial in
$L$ with coefficients - operators of lower order. The second result
is a proof that the point $\infty$ is a regular singular point
for the operator $L$ (\thref{1.2}). In the last subsection we give
an estimate of the degree $n$ (of the string equation) in terms of
the roots of the indicial equation at $\infty$. All results
 will be used in performing Darboux transformations.
Now we will fix the situation in which we are going to work.

\bde{11.1}
By $\O$ we denote the set of all functions
that are holomorphic at~$\infty$. If we write a function
$V(x)$ from $\O$  as:
\ben
   V(x) = x^{-\nu} ( a_0 + a_1\frac{1}{x} + a_2\frac{1}{x^2}
                        + \cdots ),
   a_0 \not = 0,\nu\geq 0 .
\een
then the number
\ben
   \ord (V):=-\nu
\een
will be called order of $V$ at $\infty $ and will be denoted by
$\ord (V)$.
\ede

\noindent
We introduce also the ring ${\O}[\p]$ of all
differential operators with coefficients from $\O$.
Obviously the bispectral operator $L$ is from ${\O}[\p]$.
The only properties of $L$ and $Q$ relevant for our
purposes in the present section are summed up in terms of the wave
operator as follows:

\bde{se}
We say that the operator $L\in \O [\p]$ solves the string equation
iff the following conditions are satisfied:

{\bf 1.} There is a wave operator $K=1+\a_1\p^{-1}+\cdots$
with coefficients from $\cal O$ for
which $LK=K\p^N$ and $\ord (\a_i)\leq -1$.

{\bf 2.} There is an integer $n\geq 0$ such that
$
Q = K x\p_x^{nN+1}K^{-1}
$
is a differential operator.
\noindent
 We will call the pair $(L,Q)$ a string pair. The minimal number
$n$ in {\bf 2.} will be called the string number of $L$.
\ede


\subsection{Q as a polynomial in L}


   For convenience denote by $\R$ the differential
extension of $\Cset[x]$ by adjoining the differential
indeterminates $y_1,y_2,\cdots $,(see \cite{Ka} for details). We
endow the differential ring $\R$ with graduation which will be
useful in the sequel :
for a monomial $\tau = x^{n_0}y_{i_1}^{(n_1)}\cdots
y_{i_s}^{(n_s)}$  set $\wt (\tau)=n_0-(n_1+i_1)-
\cdots -(n_s+i_s)$.
This weight provides $\R$ with the structure of a
$\Zset$-graded ring:
\ben
\R &=&\bigoplus_{n\in \Zset}{\R}_n,
\een
where ${\R}_n$ is spanned over $\Cset$ by all monomials
$\tau\in\R$ for which $\wt (\tau)=n.$ \\
This graduation can be extended in a natural way to graduation
of the ring  of all pseudo-differential
operators with coefficients from the ring $\R$, by prescribing to
the symbol of differentiation $\p$ weight $\wt (\p)=-1$. For convenience the
last mentioned ring will be denoted by $\psd \, \R$.
In this way a
pseudo-differential operator
\ben
   P=\sum_{j\leq m}a_j\p^j \   ,a_j\in \R
\een
is homogeneous of weight $n$ if for every $j$ the coefficient $a_j$
is a homogeneous element from ${\R}_{n+j}$. If $P$ is homogeneous
then by $\wt (P)$ we will denote it's weight. We will need
two lemmas. The proof of the first one
being trivial will be omitted.

\ble{11.1}
{\it (i)}
Assume that $P=\p^N + \cdots$ is a homogenous pseudo-differential
operator from $\psd \, \R$. Then  $P$
is invertible in $\psd \, \R$
and $P^{-1}$ is homogeneous with weight  $-N$. \\[5pt]
{\it (ii)} For every two homogeneous operators
$P_1$ and $P_2$ from $\psd \, \R$ with weights
respectively $n_1$ and $n_2$ their product $P_1. P_2$ is also
homogeneous and its weight is $n_1+n_2$.
\qed
\ele

\noindent
\ble{11.2}
Assume that $L_y=\p^N+\cdots $and $Q_y$ are arbitrary homogeneous
pseudo-differential operators from $\psd \, \R$. Then
one can find an integer number $n$ and homogeneous differential
operators $\tilde q_0,\tilde q_1,\ldots$ from ${\R}[\p]$ of
orders $\leq N-1$ such that:
\beqa
   Q_y  =   \tilde q_0  L_y^n  +
            \tilde q_1 L_y^{n-1} + \cdots                \label{11.2}
\eeqa
More precisely, for any $i=0,1,\ldots$, such that
$\tilde q_i \not =0 $ the weight of $\tilde q_i$
is: $\wt (Q_y) - (n-i)\wt (L_y)$.
\ele

\proof
For given $Q_y$ we will show that $\tilde q_0$ and $n$
can be determined uniquely and that they satisfy
the properties stated in the lemma.
After $\tilde q_0$ is determined we move the term
$\tilde q_0  L^n$ to the left-hand side of \eqref{11.2} and then
in the same manner we can determine $\tilde q_1$.
Now it is clear that all $\tilde q_i$ can be found successively
and the lemma will be proved .

Denote by $m$ the order of $Q_y$ and divide $m$ by $N\,:m=nN+r,
0\leq r\leq N-1$. Multiply both sides of \eqref{11.2} by $L^{-n}$
and compare the differential parts of the two pseudo-differential
operators:
\ben
      \tilde q_0=\Big( Q_yL^{-n}\Big)_+
\een
Combining this equality with \leref{11.1} we obtain the statement
of the lemma.
\qed

Denote by $\psd \, \O$ the ring of all
pseudo-differential operators with coefficients from $\O$.
To use the result of \leref{11.2} we need a ring
homomorphism
\ben
\pi : \psd \, \R \longrightarrow
      \psd \, \O
\een
defined as follows: take the unique differential homomorphism
between $\R$ and $\O$ that maps $y_i$ into $\a_i ,
i=1,2,\ldots ,$where $\a_i$ are the coefficients of the
wave operator $K$ and then extend this homomorphism to
homomorphism between $\psd \, \R$ and
                     $\psd \, \O$
by leaving $\p$ fixed. Now from the representation in
\leref{11.2} we can derive a similar one for the operators
$L$ and~$Q$.

\ble{11.3}
Let $L$ and $Q$ form a string pair and $n$
is the corresponding string number. Then
one can find differential polynomials
$\tilde q_0,\tilde q_1,\ldots ,\tilde q_n$ from $\psd \, \R$,
such that if we set $q_i=\pi (\tilde q_i) $ then:
\beqa
    Q =q_0L^n + q_1L^{n-1} + \cdots + q_n.             \label{11.3}
\eeqa
The operators $\tilde q_i$ are homogeneous. More precisely:
   $\tilde q_0=x\p_x  $ and if $\tilde q_i \not =0$, then
\ben
\wt (\tilde q_i)=-iN.
\een
The differential operator $q_n$ is not zero.
\ele

\proof
Introduce the pseudo-differential operator
\ben
  K_y=1+y_1\p^{-1}+\cdots \in \psd \, \R
\een
It is easy to check that $L_y=K_y\,  \p_x^N     \,K_y^{-1}$
and                      $Q_y=K_y\,x \p_x^{Nn+1}\,K_y^{-1}$
are homogeneous elements from
                         $\psd \, \R$
with weights respectively:
$\wt (L_y)=-N  $ and
$\wt (Q_y)=-nN $.
The definition of $\R$ was given in such a way that
$\pi (L_y)=L$ and
$\pi (Q_y)=Q$.
Applying \leref{11.2} with $L_y$,$Q_y$ we get:
\ben
   Q_y  =   \tilde q_0  L_y^n  + \tilde q_1 L_y^{n-1} + \cdots,
\een
where each $\tilde q_i$ is homogeneous with weight
$\wt (\tilde q_i)=\wt (Q_y)-(n+i)\wt (L_y)=-iN.$
Map both sides of the last equality by $\pi $:
\ben
   Q = \pi(\tilde q_0)    L^n      +\cdots +\pi(\tilde q_n)+
       \pi(\tilde q_{n+1} L_y^{-1} +\cdots\, ).
\een
Comparing the strictly pseudo-differential parts of
the operators at the two sides of the above equality
we see that:
\ben
       \pi(\tilde q_{n+1} L_y^{-1} + \cdots ) = 0.
\een
The inequality $q_n \not = 0$ holds because $n$ was
chosen to be the minimal number with the property
that $Kx\p^{nN+1}K^{-1}$ is a differential operator.
\qed


\subsection{$x=\infty$ is a regular singular point}


Here we give the proof of \thref{1.2}, i.e. that the infinite
point is regular for $L$.
Take the smallest $n$ for which $Q=K\,x\p_x^{Nn+1}\,
K^{-1}$ is a differential operator.
The idea is to assume that $x\,=\,\infty$ is irregular
for $L$ and then to assign weights to $x$ and $\p_x$ in such a way
that the most irregular terms of $L$ at $\infty$
have the highest weight. This weights will enable us
to associate with each differential operator from ${\O}[\p]$ a
$(\rho,\sigma)$-homogeneous polynomial in Y with
coefficients Laurent polynomials in $X$.
Following \cite{Dx} and using \eqref{1.6} we will get
contradiction.

We denote the ring of Laurent polynomials by $\L$.
The definition of $\rho$ and $\sigma$ is prompted by the theory of
irregular points (see e.g. \cite{BV}). Introduce the rational
number:
\ben
  r=\max (1,2+\frac{\ord V_2}{2},\cdots,2+\frac{\ord V_N}{N})
\een
(called principal level). It  can be expressed as
$\frac{r_1}{r_2}$, where $r_1$ and $r_2$ are relatively prime.
Well known fact is that $x=\infty$ is regular if and only if
$r=1$. Our
assumption that $\infty$ is irregular point yields $r>1$.
The integers $\rho=r_2$ and $\sigma=r_1-2r_2$
represent the weights of $x$ and $\p_x$
respectively. They satisfy the inequality:
\ben
 \rho +\sigma>0.
\een
The next definitions are modifications of corresponding
ones given by J. Dixmier \cite{Dx}. In the first definition we
endow the ring
${\O}[\p]$ (of differential operators with homomorphic at
$\infty $ coefficients) with $\Zset$-graded structure.

\bde{11.2}
Assume that
$L=V_0\p^n+V_1\p^{n-1}+\cdots +V_n$
is an arbitrary element of
$\D$.
For each term
$V(x)\p_x^i$
define its weight
\ben
v_{\rho,\sigma}(V(x)\p_x^i)=\rho{(\ord V)}+\sigma{i} .
\een
Then the number
\ben
v_{\rho,\sigma}(L)  :=
                      \max_{0\leq i \leq n}
                         v_{\rho,\sigma}(V_i\p^{n-i})
\een
will be called $(\rho,\sigma)$-order of $L$.
\ede

\noindent
The second definition associates to each differential
operator from ${\O}[\p]$ a $(\rho,\sigma)$-homogeneous
polynomial from ${\L}[Y]$.

\bde{11.3}
Assume the notation of the previous definition and denote by
$I(L)$ the set $\lbrace i\in\lbrace 0,1,\cdots,n\rbrace\vert
v_{\rho,\sigma}(V_i\p^{n-i})=v_{\rho,\sigma}(L) \rbrace$. The
polynomial $p\in{\L}[Y]$ defined as:\\
\begin{equation}
p=\sum_{i\in I}a_iX^{\ord V_i}Y^{n-i}  ,
\end{equation}
where $a_i\in\Cset$ are uniquely determined from the expansion
\ben
  V_i=a_ix^{\ord V _i}+(lower\ order\ terms),
\een
will be called polynomial associated with $L$.
\ede

\noindent
The following two lemmas are also taken from \cite{Dx}. Although
the situation there is slightly different the proofs are essentially the same.
We are going to prove only the first
one. The second can be proven in a similar way.

\ble{11.4}
Assume $L_1,L_2\in{\O}[\p]$ and $\rho+\sigma>0$. The
polynomial associated to the product $L_1L_2$
is the product of the polynomials associated with $L_1$ and $L_2$
respectively. The $(\rho,\sigma)$-order of this operator is:
$v_{\rho,\sigma}(L_1L_2)=v_{\rho,\sigma}(L_1)+v_{\rho,\sigma}
(L_2)$.
\ele

\proof
Set $\xi=\p_x$. Then for the product of two differential
operators we have:
\beqa
   L_1L_2 = \sum_{k=0}^{\infty}
             :\frac{\p^k L_1}{\p \xi ^k}
              \frac{\p^k L_2}{\p  x  ^k}:,             \label{11.4}
\eeqa
where : : is the normal ordering which always puts the
differentiation on the right. Write
$L_1 = a_0\xi ^{N_1}+\cdots +a_{N_1},
 L_2 = b_0\xi ^{N_2}+\cdots +b_{N_2}$. From the definition of : : we have that
\ben
:L_1L_2:=\sum_{0\leq i\leq N_1,0\leq j\leq N_2}
             a_ib_j\xi ^{N_1+N_2-i-j} .
\een
Each term in this sum satisfies the inequality
$
v_{\rho,\sigma}(a_ib_j\xi ^{N_1+N_2-i-j})\leq
v_{\rho,\sigma}(L_1)+v_{\rho,\sigma}(L_1)
$.
The equality is possible only when $i\in I(L_1)$
and $j\in I(L_2)$. On the other
hand the coefficient in front of the highest degree of $\xi $ in:
\ben
\sum_{i\in I(L_1),j\in I(L_2)}a_ib_j\xi ^{N_1+N_2-i-j}
\een
is $a_{i_1}b_{i_2}\not =0$, where $i_1$ and $i_2$ are the
minimal numbers from $I(L_1)$ and $I(L_2)$ respectively.
Thus this sum (which in fact is equal
to the product of the $(\rho,\sigma)-$ polynomials associated with
$L_1$ and $L_2$) is not zero .
The conclusion of this observations is that the
$(\rho,\sigma)$-polynomial associated with $:L_1L_2:$ is the
product of the polynomials associated with $L_1$ and $L_2$
and also
$v_{\rho,\sigma}(:L_1L_2:)=v_{\rho,\sigma}(L_1)+v_{\rho,\sigma}(L_
2 )$. To finish the proof it is enough to use formula
\eqref{11.4} and the
obvious fact that $v_{\rho,\sigma}(\p_\xi ^k L_1)\leq
v_{\rho,\sigma}(L_1)-k\sigma$ and
$v_{\rho,\sigma}(\p_x^k L_2)\leq v_{\rho,\sigma}(L_2)-k\rho$.
\qed

\ble{11.5}
Consider again two operators $L_1,L_2\in{\O}[\p]$ and denote
by $f_1$ , $f_2$ the polynomials associated with them and by $n_1$
and $n_2$ their $(\rho,\sigma)$-orders. If the fraction
$\frac{f_1^{n_2}}{f_2^{n_1}}$ is not a constant and
$\rho+\sigma>0$, then the polynomial associated with $[L_1,L_2]$
is:
\begin{equation}
\frac{\p f_1}{\p Y}\frac{\p f_2}{\p X}\
-\ \frac{\p f_1}{\p X}\frac{\p f_2} {\p Y}.                 \label{11.5}
\end{equation}
For the $(\rho,\sigma)$-order we have a
formula: $v_{\rho,\sigma}([L_1,L_2])=n_1+n_2-\rho-\sigma$.
\qed
\ele

\noindent
In order to apply these lemmas to the string equation \eqref{1.6}
we have to find the polynomials $f$ and $g$ associated with $L$ and $Q$ and
their $(\rho,\sigma)$-orders $v$ and $w$.
This requires few auxiliary results, stated in the following two lemmas.

\ble{11.6}
{\it (i).}
The $(\rho,\sigma)$-order of $L$ is $v=N\sigma$ and the polynomial
associated with $L$ is:
\ben
  f=Y^N+(at\ least\ one\ term) .
\een
{\it (ii).}
The $(\rho,\sigma)$-order of $Q$ is $w=(nN+1)\sigma +\rho$ and
the polynomial associated with $Q$ has the form:
\ben
 g=XYf^n + a_1(X,Y)f^{n-1} +\cdots +a_n(X,Y) .
\een
\ele

\proof {\it (i)}.
Since $v_{\rho,\sigma}(\p^N)=N\sigma$ the only thing we
have to check is that $ \ord (V_i)\rho
+(N-i)\sigma\leq N\sigma$ for $i=2,3,\ldots,N $
and that equality is reached
for at least one $i$. But this is obvious from the definition of $\rho$
and $\sigma$. \\

{\it (ii).}
The polynomial $g$ has the form:
\ben
 g(X,Y)=a_0(X,Y)f^n + a_1(X,Y)f^{n-1} +\cdots +a_n(X,Y),
\een
for some $a_i\in {\L}[Y]$. \leref{11.3} gives that
$a_i,i=1,2,\ldots ,n$ can have only negative degrees of $X$. But
then the coefficient at the highest degree of $Y$
in the polynomial $g^v$ is not a constant, while the corresponding
one in $f^w$ is 1. Thus
the fraction $\frac{f^w}{g^v}$ is not a constant.
Now from \leref{11.5} the polynomial $h$ associated with $[L,Q]$
is:
\ben
 h =
  \frac{\p f}{\p Y}\frac{\p g}{\p X}
   -\frac{\p f}{\p X}\frac{\p g}{\p Y}
                                       .
\een
and
$v_{\rho,\sigma}(h)=
 v_{\rho,\sigma}(f) +
 v_{\rho,\sigma}(g)-\rho -\sigma=v+w-\rho -\sigma$. On the
other hand the string equation \eqref{1.6} yields:
$v_{\rho,\sigma}(h)=(n+1)v$. From the last two relation
we derive the formula for $w$. To finish the proof
it is enough to notice that
$ v_{\rho,\sigma}(q_0L^n)=(Nn+1)\sigma +\rho. $\qed

\ble{11.7}Under the above notations
 $g=XYf^n. $
\ele

\proof
If $h$ is the polynomial associated with $[L,Q]$ then using
\leref{11.5} we have the following series of equalities:
\ben
 \rho X h = \rho X
(  \frac{\p f}{\p Y}\frac{\p g}{\p X}
   -\frac{\p f}{\p X}\frac{\p g}{\p Y}
)
=\frac{\p f}{\p Y}(wg-\sigma Y\p_Y g)
-(vf-\sigma Y\p_Y f)\frac{\p  g}{\p Y}   = \\[5pt]
  wg\frac{\p f}{\p Y}
- vf\frac{\p g}{\p Y}
=f^{-w+1}g^{v+1}\p_Y (\frac{f^w}{g^v}),
\een
where we have used that for a
$(\rho,\sigma)$-homogeneous polynomial $f$ of
$(\rho,\sigma)$-degree $v$ the following identity holds:

$$\rho X\p_X f+\sigma Y\p_Y f=vf.$$
Now the relation \eqref{1.6} leads to:
\beqa
  N\rho X. f^{n+1}                       =
  f^{-w+1}g^{v+1}\p_Y
  \Big(
       \frac{f^w}{g^v}
  \Big).                                               \label{11.6}
\eeqa
View $f$ and $g$ as elements in ${\cal K}[Y]$,
where ${\cal K}$ is an algebraic
extension of the field of fractions of the ring $\L$ containing
all the roots of the polynomials $f(Y)$ and $g(Y)$. Take
$\a (X)$ to be a zero of $f$ of order $\nu\geq 1$. Denote also
by $\mu$ the order of $\a (X)$ as a zero of $g$. Then comparing the
orders of the terms at the both sides in formula
\eqref{11.6} we obtain
\beqa
\nu(n+1)+(w-1)\nu-(v+1)\mu=\ord _{Y-\a(X)}\p_Y\Big(  \label{orders}
                         \frac{f^w}{g^v}\Big).
\eeqa
We will treat the following 2 cases separately:\\

\noindent
{\bf{case1.}} If $w\nu\not =v\mu$,
then the right side of \eqref{orders} is
$w\nu -v\mu -1$, hence $\mu=n\nu +1$.\\

\noindent
{\bf{case2.}} If $w\nu     =v\mu$, then using the formulas
for $v$ and $w$ we find
$
\mu=\frac{w}{v}\nu=\frac{(nN+1)\sigma +\rho}{N\sigma}\nu=
(n+\frac{\rho +\sigma}{N\sigma})\nu>n\nu.
$\\

\noindent
In both cases $\mu > n\nu$, which means that
$\frac{g}{f^n}$ is a polynomial in $Y$.
Now from \leref{11.6} $f^n$ divides the polynomial
$a_1f^{n-1}+\cdots +a_n$ whose degree in $Y$ does not exceed
$N-1 + N(n-1) <nN$. But the degree of $f^n$ is exactly $nN$
$\Rightarrow a_1f^{n-1}+\cdots +a_n=0$.
\qed

Now we are ready to give the proof of \thref{1.2}.\\

\noindent
{\bf{Proof of \thref{1.2}.}}
Put
$w=nv+\rho +\sigma$ and $g=XYf^n$
in \eqref{11.6} and after simplifications we get a
differential equation for $f$:
\ben
Y\p_Y(f^{\rho +\sigma})=(\rho +\sigma)Nf^{\rho +\sigma} .
\een
An immediate consequence
of this equation is that $f=c(X)Y^N$. But the choice of $\rho$ and
$\sigma$ was done in such a way that $f=Y^N +(at\,least\,one\,
term)$. This contradiction proves the theorem.
\qed

The regularity of $L$ imposes the following
restrictions on the coefficients of the
wave operator $K$.

\bco{11.2}
Let $L$ be an operator solving the string equation
and $K=1+\a_1\p^{-1}+\ldots$ is the wave operator defining the
corresponding string pair.
Then the order of the coefficient $\a_i,i=1,2,\ldots$
does not exceed the number $-i$, i.e.
\ben
   \a_i(x)\in \frac{1}{x^i}\O .
\een
\eco


\subsection{An estimate for n}

Here we want to estimate the number $n$ from the string equation
in terms of the roots of the indicial equation for $L$ at $\infty$.
For us it would be convenient to write the indicial
equation, using again the idea of the
weights. But in order to have an analogue of \leref{11.4} we
have slightly to change
the procedure of association polynomials to the elements
of ${\O}[\p]$. The next definition describes this process.

\bde{11.4}
Write every $L \in {\O}[\p] $ as
\ben
    L= V_0D^N+\cdots +V_{N-1}D+V_N
\een
where $D=x\p_x $ and assume also that:
\ben
  V_i=a_ix^{\nu _i}+(lower\,order\,terms)
\een
The number
\ben
     \wt (L):=\max _{0\leq i\leq N}\ord (V_i)
\een
will be called weight of $L$. The polynomial associated
with $L$
is from $\Cset [D]$ and is defined as follows:
\ben
  p(L)=\sum_{i:\ord (V_i)=\wt (L)}
            a_iD^{N-i}
\een
\ede

\noindent
In particular if the point $x=\infty$ is regular then
$p(\lambda)=0$ is explicitly the indicial equation
(see \cite{I}).

In terms of the above definition we can
give the following corollary from
\leref{11.3} and \coref{11.2}.

\bco{11.1}
Assume that $(L,Q)$ is a string pair and
$n$ is the corresponding string number.
Divide $Q$ by $L$ to derive
\ben
Q=Q_1L+q
\quad ,
\een
where $q$ is a differential operator of order
not exceeding $N-1$.
The weight of $q$ satisfies the inequality:
$
  \wt\, (q) \leq -nN.
$
\eco

\noindent
Obviously
$
   p(D)x^i=x^ip(D+i)
$
for every polynomial $p \in \Cset[D]$.
This observation is enough to make the following conclusion:

\ble{11.8}
Assume that
$
  L_1,L_2 \in \D
$
are two arbitrary differential operators
and denote by $f_1$ and $f_2$ the
polynomials associated with them. Then the polynomial
associated with the product $L_1L_2$ is:
\ben
  f_1(D+\wt (L_2))f_2(D).
\een
The weight of the product is a sum of the weights of the
two operators. \qed
\ele

\noindent
Now we will assume that $L$ is an operator
solving the string equation.
Denote by
$\lambda_1,\lambda_2,\ldots ,\lambda_N $
the roots of the indicial equation of $L$ at $\infty$.
The following very important fact, used
in performing Darboux transformations, is the content of
the next proposition.

\bpr{11.1}
Assume that $L$ is a differential operator
that solves the string equation.
Then we can find numbers $i$ and $j$ such that:
\ben
n\leq \frac{1}{N}\mid \lambda_i - \lambda_j \mid .
\een
\epr

\proof
Let $(L,Q)$ be a string pair. As in
\coref{11.1} divide $Q$ by $L$
\ben
Q=Q_1L+q
\een
Using that $Q$ satisfies the string equation \eqref{1.6}
we get
\beqa
 L q   =  L_1  L,                                       \label{11.7}
\eeqa
for some $L_1 \in {\O}[\p]$.\\
Denote by $f,g$ and $h$ the polynomials associated with
$L,q$ and $L_1$ respectively. For us the weight of $q$ will
be very important and will be denoted by $w$.

\leref{11.8} combined with \eqref{11.7} gives:
\ben
   f(D+w)g(D)=h(D-N)f(D).
\een
As a result we found that:
$
  f(\lambda_i +w)g(\lambda_i)=0,
$
for $i=1,2,\ldots,N. $
Using the inequality: $\deg g\leq N-1$ one can find $\lambda_i$
for which $
  g(\lambda_i)  \not  = 0 ,
$
hence
$
  f(\lambda_i +w)=0,
$
i. e. $\lambda_j=\lambda_i + w$, for some $\lambda_j$. Applying
\coref{11.1} we get that
$
  w\leq -nN.
$
This gives that
$$nN\leq \mid w \mid =
\mid \lambda _i - \lambda _j \mid .$$
\qed

\sectionnew{$\Zset _r$-invariantness of bispectral operators.}


Let ${\A}_L$ be the ring of all differential operators
commuting with $L$.
We want to prove that if the rank of $L$ is $r$ then $L$ is a
$\Zset_r-$invariant operator. The next lemma shows that it is
enough
to prove that $\Lambda(z,\p_z)$ is $\Zset_r$-invariant.

\ble{9.1}
If $\Lambda$ is $\Zset_r$-invariant then $L$ is also
$\Zset_r$-invariant.
\ele

\proof
It is enough to prove that the wave operator $K$ is
$\Zset_r$-invariant. Assume that $\Lambda$ is $\Zset_r$-
invariant. Then obviously $\Theta =b(\Lambda)$ is also
$\Zset_r$-invariant.

Now by induction on $i$ we will see that the term
$\a_i\p^{-i}$ is
$\Zset_r$-invariant. Compare the coefficients in front
of $\p^{-j}$ in the relation:
\ben
\theta                  (1        +\a_1    \p^{-1}+\cdots )   =
(1+\a_1\p^{-1}+\cdots ) (\Theta_0 +\Theta_1\p^{-1}+\cdots).
\een
Comparing the coefficients in front of $\p^0$ and $\p^{-1}$
one deduces that $\Theta =\Theta_0$ is $\Zset_r$-invariant and
that $\Theta_1 =0$.
Next assume that $\a_1\p^{-1},\a_2\p^{-2},\ldots,\a_i\p^{-i}$
are $\Zset_r$-invariant and compare the coefficients in
front of $\p^{-i-2}$:
\ben
  \theta\a_{i+2} & = &
                       \sum_{s=2}^{i+2}{\a_{i+2-s}
                            (\Theta_s      + {{s-i-2}\choose 1}
                             \Theta_{s-1}' +\cdots
                                           + {{s-i-2} \choose s}
                             \Theta_0^{(s)}              )}\\
                 &   & +\a_{i+2}\theta
                       +\a_{i+1}(\Theta_1 -\Theta_0')
\een
The last formula together with the fact that $\Theta$ is a
$\Zset_r$-invariant pseudo-differential operator and the
inductive assumption
 give that $\a_{i+1}\p^{-i-1}$ is $\Zset_r$-invariant. \qed

The next lemma shows that the algebra $A_L$ consists of
$\Zset_r$-invariant polynomials.

\ble{9.2}
Let $L$ be an operator of rank $r$. Then $A_L$ is a subalgebra
of $\Cset [z^r]$.
\ele

\proof
Let $P \in {\A_L}$. Put $b_1(P)= f(z) \in A_L$. From
\leref{2.32} we know that $f(z)$ is a polynomial. Also the degree
of $f(z)$ is a number divisible by $r$. Assume that
$f\not\in \Cset[z^r]$ and also that the coefficient
in front of the highest degree is $1$.
We can represent $f$ as:
\ben
f=f_0 +f_1,
\een
where $f_0\in \Cset[z^r]$ is formed from all terms of $f$ whose
degrees are divisible by $r$ and $f_1=f-f_0$. The polynomial
$f_0$ will be called the invariant part of $f$ and $f_1$ the
non-invariant part of $f$. Denote by $n_0$ and
$n_1$ the degrees of $f_0$ and $f_1$ respectively. Obviously
$n_0>n_1$ and $n_1$ is not divisible by r.
The idea is to construct new polynomial $\tilde f$
from $A$ in such a way that the difference $\tilde n_0 -
\tilde n_1$ between the degrees of the invariant and the non-invariant
part of $\tilde f$ is smaller. After finitely many
steps we will end up with a polynomial for which this
difference is negative, which will be a contradiction.

The polynomial $\tilde f$ can be constructed as follows:
let $n_0=kr$ and $N=pr$ set $\tilde f:=f^p-z^{kN}$. Denote
by $\tilde f_0$ and $\tilde f_1$ the invariant and the
non-invariant parts of $\tilde f$ and let
$\tilde n_0$ and $\tilde n_1$ be their degrees.
Write the following chain of equalities:
\ben
\tilde f=f^p-z^{Nk}=(f_0 + f_1)^p-z^{Nk}=f_0^p-z^{kN}
+{p \choose 1}f_0^{p-1}f_1 +\cdots.
\een
Since $pn_0=kN$ and $f_0$ is a polynomial in $z^r$ we can
conclude that $\tilde n_0\leq pn_0 -r$. The above expansion
together with $n_0>n_1$ gives that $\tilde n_1=n_0(p-1)+n_1$.
Now we can prove that the new difference is smaller:
\ben
\tilde n_0 -\tilde n_1\leq pn_0 -r-\tilde n_1
=pn_0-r-(p-1)n_0 -n_1=n_0-n_1-r.
\een
\qed

\noindent
{\bf{Proof of \thref{1.3}.}}
It remains to prove the $\Zset_r$-invariantness of $\Lambda$.
Write $\Lambda $ in the form:
\beqa
\Lambda (z,\p_z)=\sum_{i=0}^{r-1}
z^i\Lambda_i(z^r,z\p_z).                              \label{9.4}
\eeqa
where
\ben
\Lambda_i(z^r,z\p_z)=\sum_{j=0}^{n_i}\Lambda_{i,j}(z^r)
                                     (z^{nN+1}\p_z)^{n_i-j}
\een
All $\Lambda_{i,j}$ are Laurent polynomials and $n_i$ is chosen
in such a way that $\Lambda_{i,0}\not =0$, when $\Lambda_i\not
=0.$ We have to prove that all $\Lambda_i,i=1,2,\ldots ,r-1$ are $0$.
 Thus assume that at least one
$\Lambda_i \not =0$.
After applying the bispectral involution
$b_1$ on \eqref{9.4} we will get the following relation:
\beqa
\theta =\sum_{j=0  }^{n_0}Q^{n_0-j}\Lambda_{0,j}(L^{\frac{r}{N}})+
\cdots
+\sum_{j=0}^{n_{r-1}}Q^{n_{r-1}-j} \Lambda_{r-1,j}(L^{\frac{r}{N}
})L^{\frac{r-1}{N}}                                   \label{9.5}
\eeqa

The idea is to construct an operator from ${\A}_L$
whose image under
the bispectral involution is not from $\Cset[z^r]$. This will be
contradiction with \leref{9.2}. We split the construction of
such an operator into two cases:\\

\noindent
{\bf{case1.}} $n_0\leq \max\lbrace n_1,n_2,\ldots,n_{r-1}\rbrace.$\\
Denote by $\rho$ the maximal value of the numbers
$n_0,n_1,\ldots,n_{r-1}$ and by $I$ the set of all indeces $i$
for which $n_i=\rho$. Due to \leref{3.1}
\ben
(\ad \, L)^{\rho}(Q^\rho)=(\rho)!N^{\rho}L^{\rho (n+1)},
\een
hence one obtains the following relation:
\beqa
(\ad \, L)^{\rho}(\theta)=(\rho)!N^{\rho}L^{\rho (n+1)}
                  \sum_{i \in I}\Lambda_{i,0}(L^{\frac{r}{N}})
                                L^{\frac{i}{N}}.  \label{9.8}
\eeqa
Since the operator at the right-hand side commutes with $L$, it follows
that the differential operator at the left-hand side is
from ${\A}_L$. After applying the bispectral
involution to \eqref{9.8} we get that:
\ben
z^{\rho (n+1)N}\sum_{i \in I}  \Lambda_{i,0}(z^r)  z^i
\een
is an element from $A_L$. This element is not polynomial in
$z^r$ because the set $I$ includes at least one index
$i\in \lbrace 1,2,\ldots ,r-1\rbrace$.\\

\noindent
{\bf{case2.}}
$n_0 > \max \lbrace n_1,n_2,\ldots,n_{r-1}\rbrace. $\\
Now \eqref{9.5} can be written in the form:
\beqa
\theta -Q^{n_0} \Lambda_ {0,0} (L^{ \frac{r}{N} }) =
\sum_{j=1  }^{n_0-1  }Q^{n_0-j}\Lambda_{0,j}(L^{\frac{r}{N}})+
\sum_{j=0  }^{n_1    }Q^{n_1-j}\Lambda_{1,j}(L^{\frac{r}{N}})
L^{\frac{1}{N}}+
\cdots         +                                   \label{9.6} \\
\sum_{j=0}^{n_{r-1}}Q^{n_{r-1}-j}\Lambda_{r-1,j}
(L^{\frac{r}{N}})L^{\frac{r-1}{N}}                \nonumber
\eeqa
Applying $(\ad \, L)^{n_0}$ to  the above equality we see
(using \leref{3.1}) that it annihilates the operator
at the right-hand side. Hence the operator
\ben
(\ad \, L)^{n_0}(Q^{n_0} \Lambda_ {0,0} (L^{\frac{r}{N} }))=
(\ad \, L)^{n_0}(\theta)
\een
must be differential. Denote by ${N_1}$ the number
$n_0(n+1)$. Using again \leref{3.1}, i.e. that
$$
(\ad \, L)^{n_0}(Q^{n_0})=n_0!N^{n_0}L^{n_0(n+1)}
$$
we see that after
multiplying from the right both sides of \eqref{9.6} by
$L^{N_1}$ the operator on the left-hand side will
become differential. Denote this new operator by $P$. We
re-denote $\Lambda_{i,j}L^{N_1}$ by $\Lambda_{i,j}$ to avoid
complicated
notation. Thus the new relation has the form:
\beqa
P= &
     \sum_{j=1   }^{n_0-1}Q^{n_0-j}\Lambda_{0,j}(L^{\frac{r}{N}})+
     \sum_{j=0   }^{n_1  }Q^{n_1-j}\Lambda_{1,j}(L^{\frac{r}{N}})
      L^{\frac{1}{N}}+
      \cdots         + &                  \nonumber  \\
  & \sum_{j=0  }^{n_{r-1}}Q^{n_{r-1}-j}\Lambda_{r-1,j}
     (L^{\frac{r}{N}})L^{\frac{r-1}{N}}         \label{9.7}
\eeqa
We can repeat this procedure until $n_0$ is reduced to a
number smaller or equal to $\max {\lbrace}{n_1,n_2,\ldots,n_{r-1}}
                                 {\rbrace}$.
Then one proceeds as in {\bf case 1}.
\qed


\sectionnew{Darboux transformations}


In this section we will gradually simplify the operator $L$ by
successive applications of Darboux transformations.
Our goal is to obtain after
a finite number of steps a Bessel operator.

According to \thref{1.2} the point $x=\infty$ is a regular
singular point for the
operator $L$. Assume also that $L$ is a rank $r$
differential operator. From \thref{1.3} we know that in this
case $L$ is $\Zset_r-$invariant operator. Thus if we
represent $L$ as
\ben
L=\p^N+V_1\p^{N-1}+\cdots+V_{N-1}\p +V_N
\een
then the coefficients $V_i$ can be expanded as
\beqa
V_i=\frac{1}{x^i}\sum_{k=0}^{\infty}
     V_{i,k}x^{-rk}                           \label{D.1}
\eeqa
In what follows we need to split the set
$M=\lbrace \lambda_1,\lambda_2,\ldots ,\lambda_N
   \rbrace $
of roots of the indicial equation
 at~$\infty$ for $L$ into subsets of equivalent
modulo $\Zset$ numbers.

For an arbitrary set $M_i$ denote by $\lambda$ the number
in $M_i$ with minimal real part. The next lemma is a
version of a classical
result (see e. g. \cite{I}) and shows how one can pick an
$\Zset_r$-invariant function from $\ker \,  L$.

\ble{6.1}
If $\lambda$ is the minimal number of a set $M_i$, then
there is a  function $\phi_{\lambda}$ from $\ker \,  L$ which can be
expanded around $\infty$ as~:
\beq
\phi_{\lambda}(x)=x^{\lambda}\sum_{k=0}^{\infty}
                 c_kx^{-kr} \quad  ,c_0 = 1         \label{D.2}
\eeq
\ele

\noindent
We omit the proof as it repeats the classical one.

Given a function $\phi_{\lambda}$ we construct a first
order operator by setting
\ben
P_{\lambda}=\p_x - \frac{\phi_
{\lambda}'}{\phi_{\lambda}} .
\een
Then the operator $L$ can be factorised as
$L=Q_{\lambda}P_{\lambda}$ and after we perform the Darboux
transformation
\beqa
L       =Q_{\lambda}P_{\lambda}\longrightarrow
\tilde L=P_{\lambda}Q_{\lambda}                 \label{min.d}
\eeqa
the new operator $\tilde L$ will have the following properties:

\bpr{D.1}
Assume that the operator $L$ solves the string equation with a
$\Zset_r$-invariant wave operator
$K=1+\a_1\p^{-1}+\cdots,\ord (\a_i)\leq -i$. Then

(i) every operator which is obtained by a Darboux
transformation described  above also solves the string equation;

(ii) if $\lbrace \lambda_1,\lambda_2,
\ldots ,\lambda_N\rbrace $ are the roots of
the indicial equation at $\infty$
of $L$ and $\lambda=\lambda_{i_0} $ is the number with minimal
real part from some $M_i$ then the roots of the indicial equation
at $\infty$ of $\tilde L$ are
$
\tilde \lambda_k=\lambda_k -1
$
for $,k\not =i_0$ and
$
\tilde \lambda_{i_0}=\lambda_{i_0} +(N-1)
$.
\epr

\proof
Put
\beq
\tilde K=P_{\lambda}K\p^{-1}.                  \label{D.51}
\eeq
 Now we will check that $\tilde L, \tilde K$ also satisfy the
 conditions of the lemma. Let's check the first condition of \deref{se}.

\ben
\tilde L\tilde K=P_{\lambda}Q_{\lambda}P_{\lambda}K\p^{-1}=
                 P_{\lambda}L                     K\p^{-1}=
                 P_{\lambda}     K \p^N            \p^{-1}=
\tilde K\p^N.
\een

Further denote by $\tilde Q=P_{\lambda}QQ_{\lambda}$ and note
that the following sequence of equalities holds:
\ben
\tilde Q\tilde K=P_{\lambda}QQ_{\lambda}P_{\lambda}K\p^{-1}=
P_{\lambda}QLK\p^{-1} .
\een
Using the equalities $LK=K\p^N$ and $QK=Kx\p^{nN+1}$ the last
relations give
\ben
\tilde Q\tilde K=P_\lambda Kx\p^{nN+1}\p^{N-1}=
\tilde K \p x\p^{N(n+1)}=\tilde K x\p^{(n+1)N+1} +
                         \tilde K  \p^{(n+1)N  }.
\een
Now it is clear that
$\tilde K x\p^{(n+1)N+1}\tilde K^{-1}=\tilde Q-\tilde L ^{n+1}$
is a differential operator.

Denote by $g$ the polynomial associated with $P_\lambda$ and
by $h$ the one associated with $Q_\lambda$.
Obviously
$g(D)=D-\lambda$ and it has weight $-1$.
Then using
\leref{11.8} we get:
\ben
h(D-1)g(D)=(D-\lambda_1)(D-\lambda_2)\ldots (D-\lambda_N) \\
g(D-(N-1))h(D)=(D-\tilde\lambda_1)(D-\tilde\lambda_2)\ldots
   (D-\tilde\lambda_N)
\een
>From these equalities we get  the second assertion in the
lemma.
\qed

After this proposition we are close to our final goal.

\bpr{D.2}
Let $L$ be a bispectral operator with
coefficients satisfying \eqref{1.4}. Then by finitely many
$\Zset_r-$invariant
Darboux transformations we can transform it into a Bessel
operator.
\epr

\proof
We will perform  Darboux transformations in the following
way:
start with $L_0=L$. Choose an
index $i$, if there is any, for which the difference
between numbers in $M_i$ with maximal real part and
with minimal real part exceeds $N$.
Denote by ${\lambda}$ the number in $M_i$ with the minimal
real part, set $P_1=P_{\lambda}$ and factorise $L$ as
$L=R_1P_1$
then the Darboux transformation will be
\ben
L_0=L=R_1P_1 \longrightarrow L_1:=P_1R_1
\een
According to \prref{D.1} the sets $M_j^0:=M_j$ will be
transformed into sets $M_j^1$ for which the difference
between the numbers with maximal and minimal real parts
are the same for $i\not =j$. When $i=j$ there are two
cases:\\

\noindent
{\bf{case1.}} There is exactly one number in  $M_i$ with
minimal real part.
The differences between the numbers in $M_i$
are integer. Thus there is a well defined ordering:
$\lambda \geq \mu$  iff
$\lambda -\mu\geq 0$, in fact
$\lambda -\mu=\Re \lambda -\Re \mu$.
Having in mind this remark the elements of $M_i$ can
be ordered as
\ben
\lambda <\mu_1 \leq \ldots \leq \mu_s.
\een
Now the assumption about $M_i$ means
\ben
\lambda + (N-1) \leq \mu_s -1.
\een
Due to \prref{D.1} in the new set $M_i^1$ the
following inequalities must hold:
$\min M_i^1 \geq \lambda , \quad \max M_i^1 = \mu_s-1.$
Hence, the difference between the maximal and
the minimal number is reduced at least by 1.\\

\noindent
{\bf{case2.}} In $M_i$ there is at least two numbers
with minimal real part. Then the above Darboux
transformation
decreases the number of the roots with minimal real part
at least by 1.\\

\noindent
After finitely many Darboux transformations we obtain an
operator
$L_m$ such that if $M^m$ is the
set of roots of the indicial equation at $\infty$ and
$M_j^m$ are the corresponding subsets modulo $\Zset$
for $M^m$, then
\beqa
  \max \, M_j^m -\min \, M_j^m < N.               \label{D.6}
\eeqa
But again from \prref{D.1} it follows that there is an
operator
\ben
K_m=1+a_1\p^{-1}+\cdots,
\een
such that $L_mK_m=K_m\p^N$ and there is an integer
$n\geq 0$ for which
\beqa
Q_m=K_mx\p^{nN+1}K_m^{-1}                   \label{D.5}
\eeqa
is differential. The minimal $n$ with this
property, according to \prref{11.1} satisfies the
inequality:
\ben
n\leq \frac{1}{N}(  \max M_j^m -\min M_j^m) .
\een
Using \eqref{D.6} we see that $n$ must be zero. Put in
\eqref{D.5} $n=0$ and compare the differential
parts of the operators at both sides to conclude that
$$Q_m=x\p_x. $$
Now comparing the coefficients at $\p^j$ at both
sides of the string
equation \eqref{1.6} with $n=0$ we obtain the equation
\ben
-xV_j'+jV_j=NV_j.
\een
Integrating it, we obtain that
$$
V_j=v_j{x^{-N+j}}, \quad v_j \in \Cset.
$$
This shows that $L_m$ is a Bessel operator.
\qed\\

\noindent
{\bf{Proof of \thref{1.1}.}}
It remains to show that the chaine of the above Darboux
transformations can be replaced by one monomial.
First we represent the chain by following graph:
\ben
L_0=R_1P_1 \rightarrow
L_1=P_1R_1=R_2P_2 \rightarrow
L_2=P_2R_2=R_3P_3 \rightarrow
\ldots            \rightarrow
L_m=P_mR_m.
\een
If we set $A=R_1R_2\ldots R_m$ and $B=P_mP_{m-1}\ldots P_1$
then obviously:
\ben
L^m=L_0^m=R_1P_1R_1P_1\ldots R_1P_1=
          R_1  L_1^{m-1}        P_1=AB
\een
and for the Bessel operator $L_{\b'}:=L_m$
\beqa
L_{\b'}^m=P_mR_mP_mR_m\ldots P_mR_m=
       P_m  L_{m-1}^{m-1}       R_m=BA                \label{D.7}
\eeqa
The Darboux transformations do not change the rank of the
operator. Thus the rank of $L_{\b'}$ is $r$. If $r<N$
then according to \leref{77.1} there is a monomial Darboux
transformation which transforms $L_{\b'}$ into $L_\b$, where
$L_\b$ is some Bessel operator of order $r$ and rank $r$.
But the monomial Darboux transformations connecting Bessel operators
are transitive. Thus there is a monomial Darboux transformation
connecting $L$ and $L_\b$. The only thing that we have
to prove is that the operators $A$ and $B$ from
\eqref{D.7} have rational coefficients.To prove this
we need the following lemma.

\ble{D.7} Assume that $P\in {\O}[\p]$ is an operator
with holomorphic at $\infty$ coefficients. If $P$ divides
from the right some power $L_\b^d$ of a Bessel operator
\ben
L_\b=x^{-N}(D-\b_1)\ldots (D-\b_N),\quad D=x\p_x,
\een
then the coefficients of $P$ are rational.
\ele

\proof
Let $n$ be the order of $P$ and
\ben
\gamma=\b^d=(\b_1,\b_1+N,\ldots,\b_1+(d-1)N,\ldots ,
             \b_N,\b_N+N,\ldots,\b_N+(d-1)N).
\een
First we prove that $\ker \,  P$ has a basis of elements
$f_i,i=1,2,\ldots,n$ of the form:
\beqa
f_i=x^{\gamma_i}\sum_{j=0}^{r_i}p_{ij}(x)(\ln x)^j ,
                                p_{ir_i}\not =0,        \label{D.8}
\eeqa
where $p_{ij}$ are polynomials.
In general every $f\in \ker \,  P$ can be written as:
\beqa
f=\sum_{i=1}^s f_i                                     \label{D.9}
\eeqa
with $f_i$ having of the form given by \eqref{D.8} and
$\gamma_i -\gamma_j \not \in \Zset$ for $i \not =j$.
The analytical continuation around the infinite point
defines the monodromy map:
\ben
M_{\infty}:\ker \,  P \longrightarrow \ker \,  P.
\een
If an element $f=\sum_{i=1}^s f_i  $ as in \eqref{D.9} and
\eqref{D.8} is in $\ker \, P$, then
\ben
M_\infty (f)=\sum_{i=1}^s \exp (2\pi {\sqrt -1} \gamma_i)
                  x^{\gamma _i}
             \sum_{j=0}^{r_i}p_{ij}(x)(\ln x + 2\pi {\sqrt -1})^j
\een
is also in $\ker \,  P$.

Let $s$ be the minimal number for which there is an element
$f$ as in \eqref{D.9}, where none of the terms $f_i$ is in
$\ker \,  P$. From all such functions from $\ker \,  P$ with minimal
$s$ take one for which the number:
\ben
\min \lbrace r_i \mid i=1,2,\ldots, s \rbrace
\een
is minimal. We can assume that $r_s=
\min \lbrace r_i \mid i=1,2,\ldots, s \rbrace$.
Then in the following
element from $\ker \, P$:
\ben
f-\exp (-2\pi \gamma_s {\sqrt -1})M_{\infty}(f)=
\sum_{i=1}^s \tilde{f_i}                     =
x^{\tilde{\gamma_i}}\sum_{j=0}^{\tilde{r_i}}
\tilde{p_{ij}}(x)(\ln x)^j
\een
either the term $\tilde{f_s}$ vanishes (when $r_s=0$) or
 the number $\tilde{r_s}=r_s-1$ is less than $r_s$.
In both cases this is a contradiction with the choice of
$f$.

Having in mind the basis from \eqref{D.8} the action
of the operator $P$ can be written as (see \cite{I}):
\beqa
P\phi=\frac{\Wr (f_1,f_2,\ldots,f_n,\phi)}
           {\Wr (f_1,f_2,\ldots,f_n     )}.              \label{D.10}
\eeqa
Note that each derivative $f_i^{(k)}$ has the form
$f_i^{(k)}=x^{\gamma_i}F_{ik}(x,\ln x),$ where
$F_{ik}(X,Y)\in {\L}[Y]$ is a polynomial in $Y$
with coefficients - Laurent polynomials in $X$.Hence
formula \eqref{D.10} gives:
\ben
P\phi = \frac{x^{\gamma_1 +\gamma_2 +\cdots +\gamma_s}
              \sum_{i=0}^nF_i(x,\ln x)\p^i\phi           }
             {x^{\gamma_1 +\gamma_2 +\cdots +\gamma_s}
                          F_n(x,\ln x),                  }
\een
where $F_i\in {\L}[Y]$. Thus the coefficient $c_i$ in
front of $\p^i$ is
\ben
c_i=\frac{F_i(x,\ln x)}{F_n(x,\ln x)}.
\een
Since $c_i\in {\O}$ the monodromy map $M_\infty$ preserves
$c_i$. Hence
\ben
\frac{F_i(x,\ln x+2\pi{\sqrt -1}l)}{F_n(x,\ln x+2\pi{\sqrt -1}l)}=
\frac{F_i(x,\ln x)}{F_n(x,\ln x)}
\een
for every integer $l$ and also for every $l\in \Cset$
since the above equality is equivalent to an equality
between polynomials. Using that $x$ and $\ln x$
are algebraically independent we get:
\ben
\frac{F_i(X,Y+l)}{F_n(X,Y+l)}=
\frac{F_i(X,Y  )}{F_n(X,Y  )}
\een
which on the other hand is equivalent to
\ben
\frac{F_i(X,Y+l)}{F_i(X,Y)}=
\frac{F_n(X,Y+l)}{F_n(X,Y)}.
\een
Take the derivative with respect to $l$ and set $l=0$.
Then one sees that $F_i(X,Y)=c(X)F_n(X,Y)$. After
putting first $Y=\ln x , X=x$ and then $Y=0 , X=x$ it follows that
\ben
c_i(x)=\frac{F_i(x,\ln x)}{F_n(x,\ln x)}=c(x)=
       \frac{F_i(x,0  )}{F_n(x,0  )}
\een
is a rational function.
\qed


\sectionnew{Proof of the characterisation theorem.}


Essentially the proof of \thref{1.4} has already been performed in the
previous sections, as well as in \cite{BHY1, BHY5}. Below we sketch a plan how
to pick the pieces of the proof from these sources.\\

\noindent
{\bf{Proof of \thref{1.4}.}}
The implication $1)\to 3)$ is the content of \thref{1.1}. Next we consider
$3)\to 2)$. Here we use the  \deref{mon2} for monomial Darboux
transformations. If $L_\b$ is a Bessel
operator then one factorises $L_{\b}^m$ as
\beq
L_\b^m= QP, \label{8.11}
\eeq
where the operator $P$ acts on $\psi$ in the following way:
\beq
P= \frac{\Wr (f_1, f_2,\ldots, f_n, \psi)}
        {\Wr (f_1, f_2,\ldots, f_n      )}
                                                  \label{8.12}
\eeq
and the functions $f_1, \ldots, f_n$ have the structure prescribed in
\deref{mon2}. Having in mind the type of the kernel it is obvious that the
operator $P$ has only regular singularities. But then the same is true for
the operator $Q$ whose coefficients are computed by induction from the
\eqref{8.11}. Then the same is true for the product $PQ$.
At the end by the main result in
\cite{BHY1} the latter operator is bispectral.

The implication $2) \to 1)$ is trivial. The equivalence of 3) and 4) is
the content of \cite{BHY5}. We briefly describe it.

First, we recall the definition of $\W$, its subalgebras
$\WN$ and their bosonic representations introduced in \cite{BHY1}.
The algebra $w_{\infty}$ of the additional symmetries of the KP--hierarchy is
isomorphic to the Lie algebra of regular polynomial differential operators on
the circle
$$\D = \span \{ z^{\a} \partial_z^{\beta}| \; \a,  \beta
\in \Zset, \; \beta \geq 0 \}. $$
Its unique central extension \cite{KR} will be denoted by
$W_{1+\infty}$.
This algebra gives the action of the additional symmetries on tau-functions
(see \cite{ASvM, OS}). Denote by $c$ the central element of $W_{1+\infty}$
and by
$W(A)$  the image of $A \in \D$ under the natural embedding $\D
\hookrightarrow
W_{1+\infty}$ (as vector spaces). The algebra $W_{1+\infty}$ has a basis
$$c, \; J_k^l = W(-z^{l+k} \partial_z^l), \qquad l,k \in \Zset, \; l \geq 0.$$
The commutation relations of $W_{1+\infty}$ can be written most
conveniently in
terms of generating series \cite {KR}
\beq
\Bigl[ W(z^k e^{xD_z}), W(z^m e^{yD_z}) \Bigr] = ( e^{xm} - e^{yk})
W(z^{k+m}
e^{(x+y)D_z})+ \delta_{k,-m} \frac{e^{xm}-e^{yk}}{1-e^{x+y}}c,
\label{8.1}
\eeq
where $D_z = z\partial_z$.

>From the theory of KP-hierarchy it is well known
that each operator $L$ or its wave function
\eqref{2.2} defines or can be defined by the so called
{\it tau-function}, which is a function $\tau(t_1,\ldots ,
t_n,\ldots)$ in infinite number of variables $t_n, n=1,\dots$. We denote the
tau-functions of the Bessel operators $L_{\b}$ by $\tau_{\b}$.
In \cite{BHY2} a family of highest weight modules
$\M_\beta$ over $\W$ has been constructed,
using as a highest weight vector $\tau_\beta$. We briefly describe
them.

Introduce the subalgebra $\WN$ of
$\W$ spanned by $c$ and $J_{kN}^l,$ $l,k\in \Zset,$ $l \geq 0 $. It is a
simple
fact that $\WN$ is isomorphic to $\W$ (see \cite{KR}). Now put
\beq
\M_\beta = \span \Bigl\{ J_{k_1 N}^{l_1} \cdots J_{k_p N}^{l_p}
                         \tau_\beta
                 \Big| k_1 \leq \ldots \leq k_p < 0
                 \Bigr\}.                                 \label{8.22}
\eeq
The main result of \cite{BHY5} can be summed up as:

\bth{D-T}
If an element in a module $\M_\b$ is a tau-function then the corresponding
operator $L$ is a monomial Darboux transformation of some Bessel
operator $L_{\b'}$ (with eventually different $\b'$).
If an operator $L$ is
a monomial Darboux transformation of a Bessel operator $L_\b$ then the
corresponding tau-function belongs to the module $\M_\b$.
\eth

\noindent
Obviously the above cited theorem gives the equivalence between 3) and 4).
\qed


\renewcommand{\em}{\textrm}
\begin{small}
\renewcommand{\refname}{ {\flushleft\normalsize\bf{References}} }
    
\end{small}

\end{document}